\documentclass{gCOV2e}

\usepackage{epstopdf}
\usepackage{subfigure}
\usepackage{color}

\theoremstyle{plain}
\newtheorem{theorem}{Theorem}[section]
\newtheorem{corollary}[theorem]{Corollary}
\newtheorem{lemma}[theorem]{Lemma}

\theoremstyle{definition}

\theoremstyle{remark}
\newtheorem{remark}[theorem]{Remark}
\newtheorem{conjecture}[theorem]{Conjecture}

\def\N{\mathbb N}
\def\Z{\mathbb Z}	

\def\R{\mathbb R}
\def\C{\mathbb C}

\def\D{\mathbb D}

\begin{document}



\title{Hayman's classical conjecture on some nonlinear second order algebraic ODEs}

\author{
\name{Robert Conte$^{\text{ab}}$,
Tuen-Wai Ng$^{\text{b} *}$
and Cheng-Fa Wu$^{\text{b}}$
\newline
\mbox{\hspace{0.9truecm}}
\thanks{
*Corresponding author. E-mail: {\color{blue} ntw@maths.hku.hk}}}
\affil{
$^{\text{a}}$Centre de math\'{e}matiques et de leurs applications (UMR 8536),
\'Ecole normale sup\'{e}rieure de Cachan, 61, avenue du Pr\'{e}sident Wilson,
F--94235 Cachan Cedex,
France; $^{\text{b}}$Department of Mathematics,
The University of Hong Kong,
Pokfulam Road, Hong Kong}
\received{Received 10 December 2014}
}

\maketitle

\begin{abstract}
In this paper, we study the growth, in terms of the Nevanlinna characteristic function,  of meromorphic solutions of three types of second order nonlinear algebraic ordinary differential equations. We give  all their meromorphic solutions explicitly, and hence show that all of these ODEs satisfy the {\it classical conjecture} proposed by Hayman in 1996.
\end{abstract}


\begin{classcode}  34M05; 30D30; 30D35
\end{classcode}

 \section{Introduction}

 One important aspect  of  the studies of complex ordinary differential equations (ODEs) is to investigate the growth of their  solutions which are meromorphic on the whole complex plane $\C$.
  A well known problem in this direction is the following {\it classical conjecture} \cite[p.~344]{Chiang2003Halburd,Laine2008} proposed by Hayman in \cite{Hayman1996}.
\begin{conjecture}[Hayman] \label{Hayman's Classical Conjecture}
If $w$ is a meromorphic solution of
\begin{equation} \label{algebraic ordinary differential equation}
P(z,w,w',\cdots,w^{(n)})=0,
\end{equation}
where $P$ is a polynomial in all its arguments, then there exist $a,b,c \in \R^{+}$ such that
\begin{equation} \label{classical conjecture}
T(r, w) < a \exp_{n-1}(br^c), 0 \leq r < \infty,
\end{equation}
where $T(r, w)$ is the Nevanlinna characteristic of $w(z)$ and $\exp_l(x)$ is the $l$ times iterated exponential, i.e.,
\begin{equation*}
\exp_0(x) = x, \exp_1(x) = e^x, \exp_l(x) = \exp \{\exp_{l-1}(x)\}.
\end{equation*}

The conjecture for the case $n = 2$ was proposed by Bank in \cite{Bank1975}. When $n=1$, \eqref{classical conjecture} reduces to
\begin{equation} \label{finite order of growth}
T(r,w) < a r^c, 0 \leq r < \infty,
\end{equation}
and we say that a meromorphic function $w$ has finite order if it satisfies
\eqref{finite order of growth}. The infimum $\sigma$ of all possible numbers $c$ is called the order of $w$. For example, $\sin z, \cos z, \tan z, e^z $ and the gamma function $\Gamma (z)$ have order 1. The Weierstrass ellptic function $\wp(z)$, which satisfies $(\wp')^2 = 4(\wp-e_1)(\wp-e_2)(\wp-e_3)$, has order 2.
The function $e^{e^z}$ has infinite order but satisfies \eqref{classical conjecture} with $n=2, c=1$.
Throughout the paper, we use the standard notations and results of Nevanlinna theory \cite{Hayman1964,Laine1993} and we shall consider solutions meromorphic on $\C$.

In this paper, we consider the following three types of second order ordinary differential equations (ODEs)
\begin{eqnarray}
ww'' - {w'}^2 + P(w) &=& 0,
\label{second order-type 1}
\\
w'' + c {w'}^2 + P(w) &=& 0,
\label{second order-type 2}
\\
w'' + c w' + P(w) &=& 0,
\label{Generalized Fisher equation}
\end{eqnarray}
where $c \in \C$ and $P$ is a polynomial.
We prove (see Theorems \ref{type1}, \ref{type2}, \ref{second order-type 3}) that the {\it classical conjecture} holds for above equations
and moreover we find all their meromorphic solutions in closed form.

\end{conjecture}

\section{Existing results}
The {\it classical conjecture}  holds for any linear algebraic differential equation. In fact, if $f$
is a meromorphic  solution of $a_n (z) w^{(n)} + a_{n-1}(z) w^{(n-1)} + \cdots + a_{ 1}(z) w' + a_{0}(z) w = h(z)   $, where all the coefficients $a_0, a_1, \dots, a_n, h$ are polynomials, then $w$ is of finite order  \cite{Wittich1950,Chen1992Gao}.
One class of nonlinear ODEs  which supports the {\it classical conjecture} is the higher order Briot-Bouquet differential equation: $Q(w^{(n)}, w)= 0, n \in \N$, where $Q$ is a polynomial in two variables, as Eremenko, Liao and Ng \cite{Eremenko2009LiaoNg} proved that all their non-entire meromorphic solutions  belong to the class $W$, which consists of elliptic functions
and their successive degeneracies, i.e., elliptic functions, rational functions of one
exponential $\exp(kz), k \in \C$ and rational functions of $z$.

Some other positive results toward the {\it classical conjecture} are as follows.

\begin{theorem} (Gol'dberg \cite{Goldberg1956} or \cite[p.~223]{Bergweiler1998,Laine1993}) \label{Golberg1956}  
Suppose \eqref{algebraic ordinary differential equation} is a first order differential equation,  then all its meromorphic solutions   have finite order.

\end{theorem}

The {\it classical conjecture} for $n = 2$ remains open and only some partial results are known. Steinmetz \cite{Steinmetz1980} proved that \eqref{classical conjecture} holds for any second order differential equation $P (z, w, w', w'') = 0$ provided that $P$ is homogeneous in $w, w'$ and $w''$. He further proved that each meromorphic solution of such an ODE can be expressed in terms of entire functions of finite order.

\begin{theorem} (Steinmetz \cite[p.~248]{Steinmetz1980,Laine1993})\label{steinmetz1980}
 Suppose $P (z, w, w', w'') $ is homogeneous in $w,w'$ and $w''$. Then all meromorphic solutions of the ODE \eqref{algebraic ordinary differential equation} are of the form
 \begin{equation*}
 G(z) = \dfrac{g_1(z)}{g_2(z)} \exp {g_3(z)},
 \end{equation*}
 where $g_1,g_2 $ and $g_3$ are entire functions of finite order.

\end{theorem}

 For higher order cases, little is known. If we restrict ourselves to the study of entire solutions of \eqref{algebraic ordinary differential equation},
 by making use of Wiman-Valiron theory, Hayman \cite{Hayman1996} obtained a positive result to certain subclass of the algebraic differential equation

\begin{equation}\label{higher order algebraic ordinary differential equation}
P= \displaystyle \sum_{\lambda \in I} a_{\lambda}(z)w^{i_0} (w')^{i_1} \cdots (w^{(n)})^{i_n} = 0,
 \end{equation}
where $I$ consists of finite multi-indices of the form $\lambda = (i_0,i_1 \cdots, i_n)$, $i_k \in \N$ and $a_{\lambda}$ are polynomials in $z$. To state Hayman's result, we recall the definition of degree and weight of the ODE \eqref{higher order algebraic ordinary differential equation}.

The degree of each term in \eqref{higher order algebraic ordinary differential equation} is defined to be $|\lambda| = i_0 + i_1+ \cdots + i_n$ and the weight $||\lambda|| $ is defined by $||\lambda|| = i_0 + 2i_1 + \cdots + (n+1)i_n$.
We shall consider the terms with the highest weight among all those  with the highest degree in \eqref{higher order algebraic ordinary differential equation}. Let $\Lambda = \{\lambda | |\lambda| = \displaystyle\max_{\lambda' \in I} |\lambda'| \}$ and $\Omega$ be the subset of $\Lambda$ such that it consists of those terms with the highest weight, then we have

\begin{theorem}(Hayman \cite{Hayman1996}) \label{Hayman1967}
Suppose $\Omega$ is defined as above for ODE \eqref{higher order algebraic ordinary differential equation}. Let $d$ be  the  maximum degree of all the polynomials $a_{\lambda}(z)$ and suppose that
\begin{equation*}
\displaystyle \sum_{\lambda \in \Omega}a_{\lambda}(z) \not = 0.
\end{equation*}
Then all entire solutions of $(1)$ have finite order $\sigma \leq \max\{2d, 1+d\}$.
\end{theorem}
\begin{remark} The upper  bound in Theorem \ref{Hayman1967} is sharp but the same result cannot be extended to meromorphic solutions \cite{Hayman1996}.

\end{remark}

\begin{remark}
Theorem \ref{Hayman1967} includes Theorem \ref{Golberg1956} if we restrict ourselves to entire solutions of first order algebraic ODEs,  since the assumption in Theorem \ref{Hayman1967} always holds in this case (there is only one term with highest weight and highest degree).
However, it is not the case for second order ODEs. One example for which  Theorem \ref{Hayman1967} is inappliable is $w w '' - {w'}^2 - w w' = 0$,
 but it satisfies the {\it classical conjecture} according to Theorem \ref{steinmetz1980}.
\end{remark}

 Among  the second order differential equations which are not covered by Theorems \ref{steinmetz1980} and \ref{Hayman1967}, the simplest one \cite{Hayman1996,Chiang2003Halburd} is perhaps
 \begin{equation}  \label{Chaing2003Halburd}
ww'' - {w'}^2 = a_2 w''+a_1w' +a_0w + b,
\end{equation}
where $a_j  $ and $b  $ are rational in $z$ (or even constants) and are not all identically zero. If all the $a_j$ and $b$ are constants, Chiang and Halburd \cite{Chiang2003Halburd} proved that all   meromorphic solutions of \eqref{Chaing2003Halburd}   satisfy \eqref{classical conjecture} by explicitly giving all its meromorphic solutions which are either polynomials or   rational functions of one exponential. Their result is obtained mainly by combining  Wiman-Valiron theory \cite[Chapter 3]{Hayman1974,Laine1993}, local series analysis and reduction of order.  Using a different approach, Liao \cite{Liao2011}  obtained the same result. Recently, Halburd and Wang \cite{Halburd2014Wang} verified the {\it classical conjecture} for $a_j$ and $b$ rational.

\section{Main Results}
 We will show that the {\it classical conjecture}  is true for the ODEs \eqref{second order-type 1}-\eqref{Generalized Fisher equation}  which are not covered by any of the above results.
The main results  can now be stated as follows.

\begin{theorem} \label{type1}

Suppose the differential equation (\ref{second order-type 1})
\begin{equation}
ww'' - {w'}^2 + P(w) = 0,
\nonumber
\end{equation}
where $P(w) = \sum_{n = 0}^k a_n w^n, a_k \not = 0 $ is a polynomial   with constant coefficients, has non-constant meromorphic solutions, then we have  $  k \leq 4$  and its meromorphic solutions are characterized as follows:
\begin{itemize}
\item [1)] if $k = 0$ or $1$, then \eqref{second order-type 1} is included in \eqref{Chaing2003Halburd};
\item [2)] if $k = 2$, then $a_0 = a_1 = 0$ and the non-constant meromorphic solutions, which are actually the general solution, of \eqref{second order-type 1} are zero-free entire functions given by
\begin{eqnarray*}
w(z) = c_1 e^{-\frac{a_{2}}{2} z^2 + c_2 z}, \quad c_1, c_2 \in \C;
\end{eqnarray*}

\item [3)] if $k = 3 $ or $ 4$, then we must have $a_2 = 0$ and any meromorphic solution $w$ of \eqref{second order-type 1}  satisfies
\begin{equation*}
{w'}^2 + a_4 w^4 + 2a_3 w^3 + 2C w^2 - 2a_1 w - a_0 = 0, C \in \C,
\end{equation*}
whose general solution is meromorphic  \cite[Chapter 11]{Hille}  and given in the Appendix \ref{The second degree Briot and Bouquet equation}.

\end{itemize}

\end{theorem}

\begin{theorem} \label{type2}

Consider the differential equation \eqref{second order-type 2}
\begin{equation}
w'' + c {w'}^2 + P(w) = 0,
\nonumber
\end{equation}
where $c \in \C$ and   $P(w) = \sum_{n = 0}^k a_n w^n $ is a polynomial.
If the ODE \eqref{second order-type 2}
has non-constant meromorphic solutions, then we have  $  k \leq 4$ and its meromorphic solutions are characterized as follows:

\begin{itemize}
\item [1)] for  $c = 0$, we have $k \le 3$ and
\begin{itemize}
\item [i)]   non-entire meromorphic solutions of \eqref{second order-type 2} exist only for $k = 2$ or $3$ and they
are given in  the Appendix \ref{The second degree Briot and Bouquet equation} as ODE \eqref{second order-type 2} can then be reduced to ODE \eqref{second degree BB equation}.

\item [ii)] entire solutions of \eqref{second order-type 2} exist only for $k =0$ or $1$ and they are given by
\begin{eqnarray} \label{type 2-c=0,k=0,1}
w(z) =
\begin{cases}
 c_1 \sin \left(\sqrt{a_1} z\right)+c_2 \cos \left(\sqrt{a_1} z\right)-\frac{a_0}{a_1}  , \quad  k = 1, \\
 c_1 + c_2 z  - \dfrac{a_0}{2} z^2,   \quad  k = 0,
\end{cases}
\end{eqnarray}
where $c_1, c_2 \in \C$ are arbitrary constants.
\end{itemize}

\item [2)] for $c \not = 0$,
any meromorphic solution $w$ of \eqref{second order-type 2} satisfies
\begin{eqnarray*}
& 4 c^5 {w'}^2 + 4 a_4 c^4 w^4+\left(4 a_3 c^4-8 a_4 c^3\right) w^3+\left(4 a_2 c^4-6 a_3 c^3+12 a_4 c^2\right) w^2
\\
&+\left(4 a_1 c^4-4 a_2 c^3+6 a_3 c^2-12 a_4 c\right) w+4 a_0 c^4-2 a_1 c^3+2 a_2 c^2-3 a_3 c+6 a_4
 = 0 ,
\end{eqnarray*}
whose general solution is meromorphic and  given in the Appendix \ref{The second degree Briot and Bouquet equation}.

\end{itemize}

\end{theorem}

The estimate of order of meromorphic solutions of \eqref{second order-type 1} and \eqref{second order-type 2} comes as an immediate corollary to Theorem \ref{type1} and Theorem \ref{type2}.

\begin{corollary}
All meromorphic solutions of \eqref{second order-type 1} and \eqref{second order-type 2} are of finite order and hence  {\it classical conjecture} holds for ODEs  \eqref{second order-type 1} and \eqref{second order-type 2}.
\end{corollary}

Finally, using the results in \cite{Ablowitz1979Zeppetella,Gambier1910} and Lemma  \ref{solution of Fisher equation with degree 3},  we have
\begin{theorem} \label{second order-type 3}

The {\it classical conjecture} holds for the second order differential equation \eqref{Generalized Fisher equation}
\begin{equation}
w'' + c w' + P(w) = 0,
\nonumber
\end{equation}
where $c \in \C$ and   $P(w)$ is a polynomial   with constant coefficients in $w$ of degree $k  $.
\end{theorem}

Most results of this paper are contained in the third author's thesis \cite{WuChengfaThesis}.

\section{ Proof of main results}

Before proving our main results, let us recall two theorems  which will be used later.

\begin{theorem}\cite[p.~43]{Chuang1990Yang} \label{Borel type theorem1}
Let $f_j(z) $ and $g_j(z) (j=1, 2, \dots,n) (n \geq 2)$ be two systems of entire functions satisfying the following conditions:
\begin{itemize}
\item [1)] $\sum_{j=1}^n f_j(z) e^{g_j(z)} \equiv 0.$

\item [2)] For $1 \leq j, k \leq n, j \not = k, g_j(z) - g_k(z)$ is non-constant.

\item [3)] For $1 \leq j \leq n, 1 \leq h, k \leq n, h \not = k$,
\begin{equation*}
T(r, f_j) = o\{ T(r, e^{g_h - g_k})\}.
\end{equation*}
\end{itemize}
Then $ f_j(z) \equiv 0 \, (j = 1, 2, \dots, n)$.

\end{theorem}

\begin{theorem}\cite[p.~210]{Chuang1990Yang} \label{Growth estimate1}
Let $f$ and $g$ be two transcendental entire functions. Then
\begin{eqnarray*}
\lim_{r \rightarrow \infty} \dfrac{\log M(r, f \circ g)}{\log M(r, f)} = \infty, \quad \lim_{r \rightarrow \infty} \dfrac{T (r, f \circ g)}{  T(r, f)} = \infty, \\
\lim_{r \rightarrow \infty} \dfrac{\log M(r, f \circ g)}{\log M(r, g)} = \infty, \quad \lim_{r \rightarrow \infty} \dfrac{T (r, f \circ g)}{  T(r, g)} = \infty.
\end{eqnarray*}

\end{theorem}

\subsection{{\it Proof of Theorem \ref{type1}}}

If $k = 0$ or $1$, the ODE \eqref{second order-type 1} is just a special case of the ODE \eqref{Chaing2003Halburd} and thus we only need to consider the case  $k \geq 2$.

First of all,  one can check immediately that \eqref{second order-type 1} does not admit polynomial solutions for $k > 2$.
By Wiman-Valiron theory \cite[Chapter 3]{Hayman1974,Laine1993}, 
we can show that the ODE \eqref{second order-type 1} does not have any transcendental 
entire solution when $k >2$ as there is only one top degree term 
in \eqref{second order-type 1}. 
Now for a meromorphic solution $w$ of \eqref{second order-type 1} with a pole 
at some $z = z_0 \in \C$,
 we consider its Laurent series 
$w(z) = \sum_{n = 0}^{\infty} w_n (z - z_0)^{n + p}, p<0, p \in \Z, w_0 \not = 0 $.
If $k = 2$, we have  $w_0 = 0$, which is impossible and hence for $k = 2$ 
there is no nonentire meromorphic solution. 
For   $k > 2$, $w_0$ is determined by $ww'' - {w'}^2 + a_k w^k$. 
Comparing the terms determining $w_0$ yields $(k - 2) p = -2$.
Since $p \ \in \Z$ and $k > 2$, there are only two choices $3$ and $4$ for $k$. 
In conclusion, meromorphic solutions of the ODE \eqref{second order-type 1} exist only 
in the cases $k  = 2, 3$ or $ 4$ and they are entire for $k =2$ and nonentire 
for $k =3, 4$. We will divide them into the following two cases.

\begin{itemize}

\item [{\it Case 1:}] $k = 3 $ or $4$.

By considering the Laurent series expansion around a pole of $w$ (if it exists), a necessary condition for the existence of non-entire meromorphic solution of \eqref{second order-type 1} is $a_2 = 0$, otherwise logarithmic branch singularity appears in the solution.
Then with the integration factor $w^{- 3} w'$, the ODE \eqref{second order-type 1} becomes
\begin{equation} \label{reduced first order ODE1}
{w'}^2 + a_4 w^4 + 2a_3 w^3 + 2C w^2 - 2a_1 w - a_0 = 0, \quad C \in \C,
\end{equation}
whose general solution is meromorphic and is given in the Appendix \ref{The second degree Briot and Bouquet equation}.

\item [{\it Case 2:}] $k =2$.

We consider the following two subcases.

\begin{itemize}
\item [{\it Subcase 2a:}]

If a nonconstant entire solution $w$ of the ODE \eqref{second order-type 1} is zero-free on $\C$, then there exists a nonconstant entire function $h(z)$ such that $w(z) = e^{h(z)}$. Substituting $w(z) = e^{h(z)}$ into \eqref{second order-type 1} yields
\begin{equation*}
a_0  e^{-2h}  + a_1 e^{- h} + a_2 + h''= 0.
\end{equation*}
If $h$ is a transcendental entire function, then by Theorem \ref{Growth estimate1} and the  properties $T(r, a + f) = T(r,  f) + O(1), a \in \C$ and $T(r, f^{(n)}) = O(T(r, f)) $ for $r \in (0 , \infty)$ outside a possible exceptional set of finite linear measure, we have $T(r, a_2 + h'') = o\{T(r, e^{h} )\}$.
If $h(z) $ is a nonconstant polynomial, then  $T(r, a_2 + h'') = o\{T(r, e^{h} )\}$ holds obviously.
Therefore, according to Theorem \ref{Borel type theorem1}, we have
\begin{equation*}
a_0 = a_1 = a_2 + h'' = 0.
\end{equation*}

Thus \eqref{second order-type 1} has a meromorphic solution in this case only if $a_0 = a_1 =0$ and it can be explicitly solved with the
general solution given by $w(z) = c_1 e^{-\frac{a_{2}}{2} z^2 + c_2 z}$, where $c_1, c_2 \in \C$.

\item [{\it Subcase 2b:}] The  entire solution $w$ of the ODE \eqref{second order-type 1} has at least one zero in $\C$. Then $v  = \dfrac{1}{w}$ is a meromorphic function  with at least one pole in $\C$ and it satisfies
\begin{equation} \label{reduced equation2}
v v ''  - {v'}^2  - a_0 v^4 - a_1 v^3 - a_2 v^2 = 0.
\end{equation}

If $a_0$ and $a_1$ do not vanish simultaneously, from {\it Case 1}, 
we know that \eqref{reduced equation2} does not have any meromorphic solution 
since $a_2 \not = 0$. If $a_0 = a_1 = 0$, then {\it  Subcase 2a} implies $v(z)$ 
is an entire function which is a contradiction. Therefore in {\it Case 2}, 
all the nonconstant entire solutions of \eqref{second order-type 1} are zero-free.

\end{itemize}
\end{itemize}

\subsection{{\it Proof of Theorem \ref{type2}}}

If $c = 0$, then for $k\geq 4$,  by Wiman-Valiron theory and from the local series expansion around a pole of $w$ (if it exists), one can easily see   that
 \eqref{second order-type 2} does not have any entire and non-entire meromorphic solution, respectively. Hence we must have $ k \leq 3$.
  Since $c = 0,$ \eqref{second order-type 2} is a second order Briot-Bouquet differential equation and hence by \cite[Theorem 1]{Eremenko2009LiaoNg}  {\it all} the non-entire meromorphic solutions of \eqref{second order-type 2} belong to class $W$. If $P(w)$ is linear in $w$ or vanishes identically, then the general solution of the ODE \eqref{second order-type 2} is given by \eqref{type 2-c=0,k=0,1}. Suppose  $k = 2$ or $3$, then the equation \eqref{second order-type 2} has neither transcendental entire solutions (by Wiman-Valiron theory) nor polynomial solutions.  For the non-entire meromorphic solutions,
   they are given in the Appendix \ref{The second degree Briot and Bouquet equation} as ODE \eqref{second order-type 2} can be reduced to ODE \eqref{second degree BB equation}.

In the following, we shall consider the case $c \not = 0$. Similarly, we have $k \leq 4$ otherwise \eqref{second order-type 2} has no meromorphic solutions.

Assume that $ z_0 $ is neither a pole nor a  critical point of the meromorphic solution $w$ of \eqref{second order-type 2}, i.e., $w(z_0) \not = \infty, w'(z_0) \not = 0$, then there exist a neighborhood $\mathcal{N'}$ of $w_0$ and a neighborhood  $\mathcal{N}$ of  $z_0$ such that $w: \mathcal{N} \rightarrow \mathcal{N'}$ is univalent. Since $w$ is a nonconstant univalent function from  $\mathcal{N}$ to $\mathcal{N}'$, it has an inverse univalent function $z = \varphi(w)$. We define $y: \mathcal{N'} \rightarrow \C$ to be
\begin{equation} \label{transformation1}
y(w) : = w'(\varphi(w)).
\end{equation}
Therefore $y(w)$ is an analytic function in $\mathcal{N}'$.
By using \eqref{transformation1}, the ODE   \eqref{second order-type 2}  reduces to
\begin{equation} \label{transformed ODE2}
\left(\dfrac{y^2}{2} \right)' + c y^2 + P(w) = 0,
\frac{}{}\end{equation}
where $y(w)$ is defined by \eqref{transformation1} and is analytic in $w$ in some domain of $\C$.
Solving the linear ODE \eqref{transformed ODE2} yields
\begin{eqnarray}
& 4 a_4 c^4 w^4+\left(4 a_3 c^4-8 a_4 c^3\right) w^3+\left(4 a_2 c^4-6 a_3 c^3+12 a_4 c^2\right) w^2  \nonumber
\\
&+\left(4 a_1 c^4-4 a_2 c^3+6 a_3 c^2-12 a_4 c\right) w  \nonumber
\\
&+4 a_0 c^4-2 a_1 c^3+2 a_2 c^2-3 a_3 c+6 a_4+4 c^5 {w'}^2
 =  C e^{-2 c w}, C \in \C. \label{reduced first order ODE2}
\end{eqnarray}

Notice that  for nonzero $C  $, the ODE \eqref{reduced first order ODE2} has neither  non-constant polynomial solution nor transcendental entire solution (by Theorem \ref{Growth estimate1}). Next, assume $w$ is a meromorphic solution of \eqref{reduced first order ODE2} with a pole at $z = z_1$,
then $z_1$ is either a pole or a removable singularity of the l.h.s of \eqref{reduced first order ODE2} while $z_1$ is an essential singularity of $C e^{-2 c w}$ for nonzero $C$. Hence, in order for \eqref{reduced first order ODE2} to hold for $w$ meromorphic, we must have $C = 0$ for which the general solution of the equation \eqref{reduced first order ODE2} is meromorphic and given in the Appendix \ref{The second degree Briot and Bouquet equation}.
\begin{remark}
One may also apply the integration factor $e^{2 c w} w'$ for \eqref{second order-type 2} to obtain the equation \eqref{reduced first order ODE2}.
\end{remark}

\subsection{{\it Proof of Theorem \ref{second order-type 3}}}

To prove Theorem \ref{second order-type 3}, let us recall some lemmas that we will need.

\begin{lemma} \label{solution of Fisher equation with degree 2}
The equation
 \begin{eqnarray} \label{Fisher equation}
w''(z)  + c w'(z)-\frac{6}{\lambda } \left(w(z)-e_1\right) \left(w(z)-e_2\right)= 0, \lambda \not = 0
\end{eqnarray}
 has meromorphic solutions if and only if $c  \left(c^2 \lambda +25 e_1-25 e_2\right) \left(c^2 \lambda -25 e_1+25 e_2\right) = 0$ and they are given respectively as follows
\end{lemma}
\begin{itemize}
\item [(1)] if $c = 0$, then the general solution  to the equation \eqref{Fisher equation} is meromorphic and given in the Appendix \ref{The second degree Briot and Bouquet equation} as ODE \eqref{Fisher equation} can be
reduced to ODE \eqref{second degree BB equation}.

\item [(2)] for $c^2 \lambda  = 25 (e_i - e_j) \not = 0, i, j \in \{1, 2\}$, then the general solution to the equation \eqref{Fisher equation} \cite{Ablowitz1979Zeppetella,Gambier1910} is
    \begin{eqnarray} \label{solution of Fisher equation with degree 2-2}
w_2(z) &=& (e_i - e_j) e^{  \frac{- 2 c }{5} z} \wp\left( e^{ \frac{- c   }{5} z  } - \zeta_0; 0, g_3 \right) + e_j,
\end{eqnarray}
where $\zeta_0, g_3 \in \C$ are arbitrary.

\end{itemize}

\begin{lemma} \label{solution of Fisher equation with degree 3}
The ODE
\begin{eqnarray} \label{Generalized Fisher-deg 3}
w'' + c w'-\frac{2 }{\lambda ^2}\left(w-q_1\right) \left(w-q_2\right) \left(w-q_3\right)= 0, \lambda (\not = 0), c, q_1, q_2, q_3 \in \C
\end{eqnarray}
has  nonconstant meromorphic solutions if and only if $c$ satisfies
\begin{eqnarray} \label{Condition on c-gengealized Fisher-deg3}
c \prod (c \lambda + q_i + q_j - 2q_k) (- c \lambda + q_i + q_j - 2q_k) = 0,
\end{eqnarray}
where $(i j k) $ is any permutation of $(1 2 3)$, and we further have
\begin{itemize}

\item [1)] for $c = 0$, the general solution of \eqref{Generalized Fisher-deg 3}   is meromorphic and given in the Appendix \ref{The second degree Briot and Bouquet equation} as ODE \eqref{Generalized Fisher-deg 3} can be
reduced to ODE \eqref{second degree BB equation}.

\item [2)] for $c \not = 0$ satisfying \eqref{Condition on c-gengealized Fisher-deg3}, we can  classify   the nonconstant meromorphic solutions of \eqref{Generalized Fisher-deg 3} into two families
\begin{itemize}

\item [i)] for $c = \dfrac{2 q_i - q_j - q_k}{  \lambda } = \dfrac{- q_i + 2 q_j - q_k}{  - \lambda } $,
\begin{eqnarray} \label{general solutions of KPP}
w_6(z) = q_k - \frac{  q_i - q_k  } {2} e^{- \frac{  q_i - q_k }{\lambda} z}  \dfrac{\wp'(e^{- \frac{ q_i - q_k }{\lambda} z}  - \zeta_0; g_2, 0)}{\wp(e^{ - \frac{ q_i - q_k }{\lambda} z}  - \zeta_0; g_2, 0)} , \, \, \zeta_0, g_2 \, \text{arbitrary}.
\end{eqnarray}

\item [ii)] if $c = \dfrac{2 q_i - q_j - q_k}{\pm \lambda }   $,

\begin{eqnarray} \label{SP solutions of generalized fisher}
w_7(z) = \dfrac{q_j e^{\frac{q_j \left(z-z_0\right)}{\pm \lambda }}-q_k e^{\frac{q_k \left(z-z_0\right)}{ \pm \lambda }}}{e^{\frac{q_j \left(z-z_0\right)}{ \pm \lambda }}-e^{\frac{q_k \left(z-z_0\right)}{ \pm \lambda }}}, \, \, z_0 \, \text{   arbitrary},
\end{eqnarray}

which  for $q_j = q_k$ degenerates to
\begin{eqnarray} \label{rational solutions of generalized fisher}
w_8(z) = \dfrac{\pm \lambda}{z - z_0} + q_j, \, \, z_0 \, \text{   arbitrary}.
\end{eqnarray}

\end{itemize}

\end{itemize}

\begin{remark}
For $c\not = 0$,  all the meromorphic solutions of the equation \eqref{Generalized Fisher-deg 3} are given by  \eqref{general solutions of KPP}-\eqref{rational solutions of generalized fisher}  and the solution \eqref{general solutions of KPP} is the general solution.
\end{remark}

\end{lemma}

\noindent {\it Proof   of Lemma \ref{solution of Fisher equation with degree 3}.}
One can easily see that constant solutions of the ODE \eqref{Generalized Fisher-deg 3} are $w = q_n, n=1,2,3  $. Next we consider  nonconstant meromorphic solutions of \eqref{Generalized Fisher-deg 3}.
By making use of Wiman-Valiron theory, it can be proven immediately that the ODE \eqref{Generalized Fisher-deg 3} does not have any nonconstant transcendental entire solution. Meanwhile, the ODE \eqref{Generalized Fisher-deg 3} does not admit any nonconstant polynomial solution. Consequently, each nonconstant meromorphic solution of the equation \eqref{Generalized Fisher-deg 3} should have at least one pole on $\C$.

Suppose $w$ is a meromorphic solution of \eqref{Generalized Fisher-deg 3} with a pole at $z = z_0$. Without loss of generality, we may assume $z_0 = 0$ then $w(z) = \sum_{j=p}^{+\infty} w_j z^{j},  - p \in \N, w_p \not = 0$.
Substituting the series expansion of $w$ into the ODE \eqref{Generalized Fisher-deg 3} gives $p = - 1,w_{-1} =  \pm \lambda$.
The ODE \eqref{Generalized Fisher-deg 3} has Fuchs indices $-1, 4$ and the corresponding compatibility conditions regarding the existence of meromorphic solution are
\begin{eqnarray} \label{Compatibility Condition-Generalized Fishder-deg 3}
\begin{cases}
  c \prod (c \lambda + q_i + q_j - 2q_k)  = 0, \text{ if }  w_{-1} = \lambda,
\\
 c \prod   (- c \lambda + q_i + q_j - 2q_k)  = 0, \text{ if } w_{-1} = -\lambda,
\end{cases}
\end{eqnarray}
where $(i j k) $ is any permutation of $(1 2 3)$.

Now we compare the ODE \eqref{Generalized Fisher-deg 3} with the following second order ODE
\begin{eqnarray} \label{decomposition1}
[D - f_2(w)] [D - f_1(w)] (w - \alpha) = 0,
\end{eqnarray}
where $D = \dfrac{d}{d z}, \alpha \in \C$ and $f_i(w) = A_i w + B_i, A_i, B_i \in \C, i = 1, 2$.
Expanding \eqref{decomposition1} gives
\begin{eqnarray} \label{decomposition2}
w'' - (f_1 + f_2 + \dfrac{d f_1}{d w} w - \alpha \dfrac{d f_1}{d w}) w' + f_1 f_2  (w - \alpha) = 0.
\end{eqnarray}
Identifying the equations \eqref{Generalized Fisher-deg 3} and \eqref{decomposition2}   leads to the conditions
\begin{eqnarray} \label{condition-decomposition}
\begin{cases}
f_1 + f_2 + \dfrac{d f_1}{d w} w - \alpha \dfrac{d f_1}{d w}+c = 0,
\\
f_1 f_2  (w - \alpha) = - \frac{2}{\lambda ^2} \left(w -q_1\right) \left(w -q_2\right) \left(w -q_3\right).
\end{cases}
\end{eqnarray}

 One can check that the compatibility conditions \eqref{Compatibility Condition-Generalized Fishder-deg 3} hold  if and only if the conditions \eqref{condition-decomposition} are satisfied or $c = 0$.

 If $c = 0$, then the ODE \eqref{Generalized Fisher-deg 3} reduces to a first order Briot-Bouquet differential equation through multiplying it by $w'$ and integration. Therefore all its meromorphic solutions  belong to class $W$
 and they are given in the Appendix \ref{The second degree Briot and Bouquet equation}.

 For $c\not = 0$ and assuming $\eqref{Compatibility Condition-Generalized Fishder-deg 3}$ from now on, due to the symmetry in \eqref{Compatibility Condition-Generalized Fishder-deg 3} and the fact that $w$ has at least one pole on $\C$, it suffices to consider the case $c = (-q_1+2 q_2-q_3)/\lambda  \not = 0$ and one  choice for $A_i, B_i, i = 1, 2$ and $\alpha$ is
 \begin{eqnarray}  \label{parameters of the decomposition}
A_1 = -\frac{1}{\lambda },  A_2 = \frac{2}{\lambda }, B_1 = \frac{q_3}{\lambda }, B_2 = -\frac{2 q_2}{\lambda },\alpha = q_1.
\end{eqnarray}

As a consequence, if $c = (-q_1+2 q_2-q_3)/\lambda \not = 0$, then \eqref{Generalized Fisher-deg 3} can be written as
\begin{eqnarray} \label{decomposition3}
[D - \frac{2  }{\lambda }w -B_2] [D + \frac{w  }{\lambda } - B_1] (w - \alpha) = 0.
\end{eqnarray}

Let $G(z) = [D + \dfrac{w  }{\lambda } - B_1] (w - \alpha)$, then we have $[D - \dfrac{2  }{\lambda }w -B_2] G(z) = 0$ from which one can solve for $G(z) = \beta e^{\int \frac{2  }{\lambda } w dz} e^{B_2 z}  , \beta \in \C$. If $\beta = 0$, from $G(z) = [D + \dfrac{w  }{\lambda } - B_1] (w - \alpha) = 0$,
we are able to obtain the first family of meromorphic solutions of the ODE \eqref{Generalized Fisher-deg 3}
\begin{eqnarray} \label{pargicular solution-Generalized fisher}
w(z) = \dfrac{q_1 e^{\frac{q_1 \left(z-z_0\right)}{\pm \lambda }}-q_3 e^{\frac{q_3 \left(z-z_0\right)}{ \pm \lambda }}}{e^{\frac{q_1 \left(z-z_0\right)}{ \pm \lambda }}-e^{\frac{q_3 \left(z-z_0\right)}{ \pm \lambda }}}, \, \, z_0 \in \C.
\end{eqnarray}

For $\beta \not = 0$, we let $H(z) = e^{\int \frac{2  }{\lambda } u dz}$ which satisfies $H'(z) = 2 w(z) H(z) / \lambda$ and
\begin{eqnarray} \label{decomposition4}
[D + \frac{u  }{\lambda } - B_1] (w - \alpha) = \beta e^{B_2 z}   H(z),
 \end{eqnarray}
hence, $w $ is meromorphic if and only if $H$ is meromorphic. By the substitution of $w = \dfrac{\lambda}{2} \dfrac{H'}{H}$ into \eqref{decomposition4}, we have
\begin{eqnarray} \label{decomposition5}
-2 B_1 \lambda  H H'+4 \alpha  B_1 H^2-4 \beta  e^{B_2 z} H^3+2 \lambda  H H''-2 \alpha  H H'-\lambda  {H'}^2=0.
\end{eqnarray}
If we let $H(z) = e^{-B_2 z} h(z)$, then the ODE \eqref{decomposition5} reduces to
\begin{eqnarray} \label{decomposition6}
\left(2 B_1+B_2\right)  \left(2 \alpha +B_2 \lambda \right) h^2 -2 h  \left(\left(\alpha +\left(B_1+B_2\right) \lambda \right) h' -\lambda  h'' \right) -\lambda  {h'}^2-4 \beta  h ^3=0.
\end{eqnarray}

Suppose $h$ is a meromorphic solution of \eqref{decomposition6}. W.L.O.G, we assume that it has  a pole at $z = 0$ and $h(z) = \sum_{j=p}^{+\infty} h_j z^{j}, - p \in \N, h_p \not = 0$  then one can check that $p = -2$ and the Fuchs indices of the ODE \eqref{decomposition6} are $- 1, 4$ with the compatibility condition
\begin{eqnarray} \label{Com Condition-h(z)-1}
\left(\alpha +\left(B_1+B_2\right) \lambda \right){}^2 \left(2 \alpha \lambda \left(10 B_1+B_2\right)   -8 \alpha ^2 +\left(-8 B_1^2+2 B_2 B_1+B_2^2\right) \lambda ^2\right)=0,
\end{eqnarray}
which by the substitution of \eqref{parameters of the decomposition} reduces to
\begin{eqnarray} \label{Com Condition-h(z)-2}
 \left(q_1+q_2-2 q_3\right) \left(2 q_1-q_2-q_3\right) \left(q_1-2 q_2+q_3\right) = 0,
 \end{eqnarray}
which implies $q_3 =   \left(q_1+q_2\right) / 2$ or $q_1 =   \left(q_2+q_3\right) / 2$ since $c = (2 q_2 - q_1 - q_3)/\lambda \not = 0 $.

Then by the substitution of \eqref{parameters of the decomposition}, the ODE \eqref{decomposition6} reduces to
\begin{eqnarray*}
&&\begin{cases}
-\lambda ^2 h'(z)^2+\lambda  h(z) \left(2 \lambda  h''(z)+3 \left(q_2-q_1\right) h'(z)\right)+2 \left(q_2-q_1\right){}^2 h(z)^2-4 \beta  \lambda  h(z)^3=0,
\\
q_3 = \frac{1}{2} \left(q_1+q_2\right).
\end{cases}
\\
&&\begin{cases}
-\lambda ^2 h'(z)^2+\lambda  h(z) \left(2 \lambda  h''(z)+3 \left(q_2-q_3\right) h'(z)\right)+2 \left(q_2-q_3\right){}^2 h(z)^2-4 \beta  \lambda  h(z)^3= 0,
\\
q_1 = \frac{1}{2} \left(q_2+q_3\right).
\end{cases}
\end{eqnarray*}

Next, it suffices to consider the case $q_3 =  \left(q_1+q_2\right) / 2$ due to the symmetry in the above two equations. By the translation against the dependent variable $u$, we may further assume $q_3 = 0$ which implies $q_1 + q_2 = 0$. Let us come back to equation \eqref{decomposition5}, which by the substitution of \eqref{parameters of the decomposition} with $q_1 = -q_2 \not = 0, q_3 = 0$ reduces to
\begin{eqnarray} \label{reduction1 of generalized fisher}
-\lambda  H'(z)^2+2 H(z) \left(\lambda  H''(z)+q_2 H'(z)\right)-4 \beta  H(z)^3 e^{-\frac{2 q_2 z}{\lambda }} = 0.
\end{eqnarray}
Performing the transformation $H (z) = v(\zeta), \zeta = e^{-\frac{q_2}{\lambda }z } $ gives
\begin{eqnarray} \label{reduction5 of generalized fisher}
{v'}^2 - \dfrac{2 \beta  \lambda} {q_2^2}v^3 + C v  = 0,
\end{eqnarray}
whose general solution is meromorphic and can be found in the Appendix \ref{The second degree Briot and Bouquet equation}.

Finally, for  $c = (- q_1 + 2 q_2  - q_3)/\lambda  \not = 0$ and $ q_3 =  \left(q_1+q_2\right)/2$,
which implies $c = - (2 q_1 - q_2 - q_3) / \lambda$,  we obtain the meromorphic solutions of the ODE \eqref{Generalized Fisher-deg 3}
(which meanwhile is the general solution)
\begin{eqnarray} \label{general solutions of generalized fisher}
w(z) = - \dfrac{  q_2 - q_3  } {2} e^{- \frac{ q_2 - q_3 }{\lambda} z}  \dfrac{\wp'(e^{- \frac{ q_2 - q_3 }{\lambda} z}  - \zeta_0; g_2, 0)}{\wp(e^{ - \frac{q_2 - q_3 }{\lambda} z}  - \zeta_0; g_2, 0)} + q_3, \, \, \zeta_0, g_2 \in \C.
\end{eqnarray}

\begin{lemma}  (\cite[p.~5]{Laine1993}) \label{growth comparsion--linear measure}
Let $g:(0,+\infty)\rightarrow \R$ and $h:(0,+\infty)\rightarrow \R$ be monotone increasing functions such that $g(r) \leq h(r)$ outside of an exceptional set $F$ with finite linear measure. Then, for any $\alpha > 1$, there exists $r_0 > 0$ such that $g(r) < h(\alpha r )$ holds for all $r \geq r_0$.
\end{lemma}

\noindent {\it Proof} \, of Theorem \ref{second order-type 3}.
From the expression of solutions in Lemmas \ref{solution of Fisher equation with degree 2} and \ref{solution of Fisher equation with degree 3}, it suffices to focus on $w_2(z)$ and $w_6(z)$ because other solutions belong to the class $W$ which only consists of meromorphic functions of finite order.

We claim that for every $\alpha \, (\not =0) \in \C $, there exists $A, B \in \R^{+}$ such that
\begin{equation*}
T(r, \wp(e^{\alpha z}; \omega_1, \omega_2) ) <  A \exp(B r ), 0 \leq r < \infty,
 \end{equation*}
 where $\omega_1, \omega_2 \in \C \backslash \{0\} \, (\omega_1/\omega_2 \not \in \R)$ are the periods of $\wp(z)$. Since $\wp (z ; \omega_1, \omega_2) =  \wp (z/\omega_1; 1, \tau) / \omega_1^2$, where $\tau =  \omega_2/\omega_1$, we only need to prove the claim for $\wp(e^{\alpha z}; 1, \tau)$, where $ \tau \in \mathbb{H}$. For brevity, we  denote $\wp(e^{\alpha z}; 1, \tau)$ by $\wp(e^{\alpha z})$.

From the theory of elliptic functions \cite{MagnusSoni1966Oberhettinger,Segal1981}, we know that $\wp( z; 1, \tau)$ satisfies the first order ODE 
${\wp'}^2 = 4 (\wp - e_1)  (\wp - e_2)  (\wp - e_3)$, where $ e_1 = \wp(1/2), e_2 = \wp(\tau/2), e_3 = \wp((1+ \tau)/ 2) $ are distinct.
Next we consider $\overline{N}_{\wp(e^{\alpha z})}(r, e_j) : = \displaystyle{\int_0^{r} \frac { \overline{n}_{\wp(e^{\alpha z})}(t, e_j)}{t} dt}, j = 1, 2, 3$,  where $\overline{n}_f(r,a) $
denotes the number of poles of  $ 1/( f - a) $ in  $ \D(r) = \{ z \in \C | |z|< r \}$, without counting multiplicity.

Let $T = 2 \pi i/\alpha = |T| e^{i \beta}, \beta \in [0, 2 \pi)$, $R$ be the region enclosed by the rectangle $\{z \in \C | z = x + i y, 0 \leq |x|, |y| <r \}$ and $R' = e^{-i \beta} R = \{z' \in \C| z' =  e^{-i \beta} z, z \in R\} \supset \D(r)$.

Then we have $\overline{n}_{\wp(e^{\alpha z})}(t, e_j) \leq \overline{n'}_{\wp(e^{\alpha z})}(t, e_j) $, where $\overline{n'}_f (r,a) $
denotes the number of poles of  $ 1/( f - a) $ in  $R'$, without counting multiplicity.
As $e^{\alpha z}$ has a period $T$, we have $\overline{n'}_{\wp(e^{\alpha z})}(t, e_j)  \leq \dfrac{2 ([t] + 1)}{T} \overline{n''}_{\wp(e^{\alpha z})}(t, e_j) $, where $[r]$ is the integer part of $r \geq 0$ and $\overline{n''}_{f}(t, a) $ is the number of poles of  $ 1/( f - a) $ in  $R_1 = \{e^{- i \beta }z | z= x + iy,   -t < x < t, 0 \leq y < T \}$, without counting multiplicity.
Since $\wp^{-1}(e_1) = \{\frac{1}{2} + m + n \tau | m, n \in \Z\}$, $\overline{n''}_{\wp(e^{\alpha z})}(t, e_1) \leq (2[e^{|\alpha| t}] +1) \times \dfrac{2e^{|\alpha| t}} {|\tau|}$. Therefore,
\begin{eqnarray*}
\overline{N}_{\wp(e^{\alpha z})}(r, e_1) &=& \displaystyle{\int_0^{r} \frac { \overline{n}_{\wp(e^{\alpha z})}(t, e_1)}{t} dt}
\\
&=& \displaystyle{\int_{\delta}^{r} \frac { \overline{n}_{\wp(e^{\alpha z})}(t, e_1)}{t} dt}
\\
&\leq& \displaystyle{\int_{\delta}^{r} \frac { 2 ( [t] + 1 )}{T t} \overline{n''}_{\wp(e^{\alpha z})}(t, e_1)dt}
\\
&\leq& \displaystyle{ \frac { 4 ( 1 + \frac{1}{\delta} )}{T |\tau|} \int_{\delta}^{r}e^{|\alpha| t} (2 e^{|\alpha| t} +1)dt}
\\
&\leq& \displaystyle{ \frac { 4 ( 1 + \frac{1}{\delta} )}{T |\tau|} \dfrac{e^{2 |\alpha| r} + e^{|\alpha| r} - 2 }{|\alpha|}}
\\
&<& a_1 e^{b_1 r},
\end{eqnarray*}
where $a_1 =  \dfrac { 4 ( 1 + \frac{1}{\delta} )}{T |\alpha \tau|} > 0, b_1 = 2 |\alpha| > 0$, and $\delta > 0$ is chosen such that $\wp(e^{\alpha z})$ omits $e_1$ in $\D(\delta)$. Applying the same argument, we can obtain the upper bounds of $ \overline{N}_{\wp(e^{\alpha z})}(r, e_j), j =2 ,3$ which are given by
\begin{eqnarray*}
\overline{N}_{\wp(e^{\alpha z})}(r, e_j)  < a_j e^{b_j r}, 0< a_j, b_j ,  j = 2, 3.
\end{eqnarray*}
According to the Second Main Theorem of Nevanlinna theory, we have
\begin{eqnarray*}
T(r, \wp(e^{\alpha z}  )) &\leq& \sum_{j=1}^3\overline{N}_{\wp(e^{\alpha z})}(r, e_j) + S(r,\wp(e^{\alpha z}  )),
\end{eqnarray*}
where $S(r,\wp(e^{\alpha z}  )) = o(T(r, \wp(e^{\alpha z}  ))$, for all $r \in [0,+ \infty) $   outside an exceptional set $E \subset (0,+ \infty) $ with finite linear measure. Hence,
\begin{eqnarray*}
T(r, \wp(e^{\alpha z}  )) < a' e^{b'r}
\end{eqnarray*}
holds for all $r \in [0,+ \infty) - E$, where $a' = (1 + \varepsilon) (a_1 + a_2 + a_3), \varepsilon >0, b' = \displaystyle{ \max_{1\leq j \leq 3}{b_j}}$. According to Lemma \ref{growth comparsion--linear measure}, for $\gamma = 3/2$, there exists $r_0 > 0$ such that $T(r, \wp(e^{\alpha z}  )) < a' e^{3b'r/2}$ for all $r \geq r_0$.  On the other hand, it is obvious that there exist $a'', b'' >0$ such that $T(r, \wp(e^{\alpha z}  )) < a'' e^{b''r}$ for all $0 \leq r < r_0$. As a consequence, we have
\begin{eqnarray*}
T(r, \wp(e^{\alpha z}  )) < A e^{B r}, 0 \leq r < \infty,
\end{eqnarray*}
where $A = \max\{a', a''\}, B = \max\{3b'/2,b''\}$.

From the proof of our claim, it is easy to see that the same conclusion holds for $\wp(k_1 \exp\{\alpha z\} + k_2), \wp'(k_1 \exp\{\alpha z\} + k_2), k_1, k_2, \alpha \in \C$ as well. By making use of the properties $T(r, f g) \leq T(r,f) + T(r,g)$ and $T(r, f) = T(r, 1 / f) + O(1)$, we conclude that there exists $a, b>0$ and $ c = 1$, such that
\begin{eqnarray*}
T(r, w_i(z)) < a e^{br}, 0 \leq r < \infty, i = 2, 6.
\end{eqnarray*}
Thus,  the proof is complete.

\section*{Funding}
The first author was partially supported by PROCORE -- France/Hong Kong joint research scheme under [F--HK39/11T, HKU 704409P]; RGC grant under [HKU 704611P].
The second author was partially supported by PROCORE -- France/Hong Kong joint research scheme under [F--HK39/11T]; RGC grant under [HKU 704409P].
The third author was partially supported by RGC grant under [HKU 704409P, HKU 704611P]; a post-graduate studentship at HKU.

\appendices

\section{The second degree Briot and Bouquet equation} \label{The second degree Briot and Bouquet equation}

We recall here the various expressions for the general solution of the first order second degree
binomial equation of Briot and Bouquet
\begin{eqnarray} \label{second degree BB equation}
& &
\left(\frac{d u}{d z}\right)^2=a_k \prod_{j=1}^k (u-e_j), \ 
\end{eqnarray}
in which the integer $k$ runs from $0$ to $4$,
$a_k$ is a nonzero complex constant and $e_j$ are complex constants.
Let us denote $\wp$ and $\zeta$ the functions of Weierstrass,
\begin{eqnarray*}
& &
{\wp'}^2=4 \wp^3-g_2 \wp - g_3,\ \zeta'=-\wp.
\end{eqnarray*}
The general solution is (the arbitrary origin $z_0$ of $z$ is omitted) \cite[Table 1, p.~73]{Murphy1960},
\begin{eqnarray}
& &
\left\lbrace
\begin{array}{ll}
\displaystyle{
u=a_4^{-1/2}\left(\zeta(z+a)-\zeta(z-a)-2\zeta(a)+A\right) = a_4^{-1/2}\left(A - \dfrac{\wp'(a)}{\wp(z) - \wp(a)} \right),
}\\
\displaystyle{ \qquad k=4,\ e_j \hbox{ all distinct},}
\\
\displaystyle{
\frac{1}{u-e_1}=A \cosh(B z)+C,\                           k=4,\ \hbox{ one double root } e_1,
}\\ \displaystyle{
u=\frac{e_1+e_2}{2}+A \coth(B z),\                         k=4,\ \hbox{ two double roots } e_1,e_2
}\\ \displaystyle{
\frac{u-e_4}{u-e_1}=A\left[\frac{e_4-e_1}{2} z \right]^2,\ k=4,\ \hbox{ one triple root } e_1,
}\\ \displaystyle{
u=e_1\pm\frac{a_4^{-1/2}}{z},\                               k=4,\ \hbox{ one quadruple root},
}\\ \displaystyle{
u=\frac{e_1+e_2+e_3}{3}+\frac{4}{a_3} \wp(z,g_2,g_3),\     k=3,\ e_j \hbox{ all different},
}\\ \displaystyle{
u=A+B \coth^2(C z),\                                       k=3,\ e_1=e_2\not=e_3,
}\\ \displaystyle{
u=e_1+\frac{4}{a_3 z^2},\                                  k=3,\ e_1=e_2=e_3,
}\\ \displaystyle{
u=\frac{e_1+e_2}{2}+A \cosh(B z),\                         k=2,\ e_1\not=e_2,
}\\ \displaystyle{
u=e_1+e^{\pm \sqrt{a_2} z},\                               k=2,\ e_1=e_2,
}\\ \displaystyle{
u=e_1+\frac{a_1}{4}z^2,\                                   k=1,\
}\\ \displaystyle{
u=\pm\sqrt{a_0}z,\                                         k=0,\
}
\end{array}
\right.
\label{eqBB21111Solution}
\end{eqnarray}
in which $\wp(a),A,B,C$ are algebraic functions of $a_k,e_j$.

\end{document}